\documentclass[11pt]{article}
\usepackage[T1]{fontenc}
\usepackage[font=small]{caption}
\usepackage[margin=1in]{geometry} 
\usepackage{amsmath,amssymb}  
\usepackage{array}
\usepackage{float}
\usepackage{booktabs}
\usepackage{graphicx}
\usepackage{xcolor}           
\usepackage{hyperref}
\usepackage{url}
\usepackage{todonotes}
 
\usepackage{authblk}

\title{Leveraging Classical and Quantum Computing for Process Systems Engineering Applications: Decomposition Algorithm with Ising Solvers for Efficient Discrete Landscape Exploration}
\date{} 
\author{Yirang Park}
\author{David E. {Bernal Neira}}
\affil{Davidson School of Chemical Engineering, Purdue University, 480 Stadium Mall Drive, West Lafayette, IN 47907, USA}

\begin{document}

\maketitle
 
\begin{abstract}
\noindent

Conceptual process design is a crucial aspect of chemical engineering that involves process synthesis in consideration of the various process phenomena. 
Mixed-integer nonlinear programming (MINLP) offers a powerful framework for modeling such design problems by combining discrete and continuous variables; however, the combinatorial complexity of discrete choices, coupled with nonlinearities, makes MINLPs challenging to solve as monolithic problems. 
Various simulation and mathematical optimization tools have been developed for solving these difficult problems, often utilizing decomposition strategies that break the original problem into more manageable subproblems.
Additionally, recent advances in computational hardware create opportunities for addressing different parts of the problem more efficiently.
Discrete subproblems can potentially benefit from quantum Ising solvers, while simulators and nonlinear solvers offer powerful tools for handling nonlinearities.
To fully exploit these emerging computational capabilities, we propose an integrated approach that decomposes MINLP problems into discrete and continuous components and solves each subproblem using the most suitable computational method.
In this work, two case studies are presented: an illustrative example involving the selection of an ionic liquid and its process design, and a more complex problem of drug substance manufacturing process optimization.
The discrete subproblem in each case is formulated as an integer programming problem and solved using a commercial classical optimization solver.
For comparative analysis, the same subproblem is reformulated as a quadratic unconstrained binary optimization and solved with various Ising solvers including simulated annealing (SA), quantum annealing (QA), and entropy computing (EC).
The continuous subproblem is solved using a classical optimizer and a simulator-based optimization approach, respectively.
In both discrete subproblems, the commercial classical solver achieved the shortest runtime in terms of computational efficiency, whereas EC took the longest, followed by QA and SA, in reaching feasible and optimal solutions.
The heuristic methods identified all or most of the feasible solutions in a single run, demonstrating advantages in solution diversity and more efficient and broad exploration of the solution space over the classical solver. 
In contrast, the global solver provides a global optimality guarantee and rapid convergence speed. 
In the context of process design, where insights of alternative feasible designs and comparisons of costs are more valuable than a single optimal solution, heuristics methods offer a better suited approach for decision-making.
This comparative analysis highlights the distinct strengths of each method and underscores the potential of this heterogeneous computing approach that leverages different methods to address practical optimization problems.  
\end{abstract}

\section{Introduction}

Mixed-Integer Nonlinear Programming (MINLP) is an optimization framework that combines both discrete (integer or binary) and continuous decision variables, incorporating nonlinear relationships in the objective function and/or constraints, thereby enabling the modeling and solution of complex, real-world problems.
Many applications, such as process design, operations research, and finance, involve decision-making problems that can be effectively modeled as MINLP problems \cite{Belotti2013,minlpreview2019}. 
Mathematical programming and MINLP formulations have been widely applied to multiple domains in chemical engineering, including flowsheet synthesis~\cite{Turkay1996,Chen2017}, heat- and utility-integration~\cite{YEE1990,CIRIC1991}, and supply chain optimizations~\cite{GrossmannEO2005,You2010}, as a means to represent the process systems for integrated decision-making.
The systematic and rigorous approach of these optimization-based methods is particularly advantageous in identifying the optimal process configuration in the  process design space, and substantial research efforts have been made spanning superstructure generation, model formulation, and algorithmic development in the field \cite{Mencarelli_2020,linan_trends_2025}.

Solving MINLP problems remains a challenge, as these problems are complex and many practical instances are non-deterministic polynomial-time (NP) hard. Finding a good or even a feasible solution can be computationally demanding \cite{Liberti2011}.
Various practical methods have been developed, and solvers have advanced to address this complexity; however, the problems and their scales that can be solved are still limited and smaller than what can be modeled \cite{minlpreview2019}.  
One of the most common approaches to solving MINLP problems is to decompose the problem into smaller subproblems, typically composed of discrete (i.e., process configuration) and continuous (i.e., process operation) parts.
The discrete subproblem can be solved to find candidate values for the integer variables, which can then be fixed, and the remaining problem can be solved as a nonlinear program (NLP) or passed to a simulator for evaluation.
This divide-and-conquer strategy enables more efficient and flexible use of specialized solvers and methods tailored to the subproblems' complexities, thereby making the overall problem more tractable.

While the nonlinear continuous subproblems are well-handled by the advanced nonlinear algorithms and solvers or simulation tools, the combinatorial complexity of the discrete subproblems remains a computational challenge \cite{constante2025}.
A discrete optimization problem can be written as  
\begin{equation}
\min_{x\in\{0,1\}} \quad f(x) \quad \text{s.t.} \quad x \in X
\end{equation}
where $X$ is a set of constraints that define the feasible space of the problem.
A branch-and-bound (BB) method is a classical, deterministic approach to solving discrete optimization problems, which systematically explores the solution space by following a search tree and eliminating infeasible branches \cite{Land1960}.
Namely, \texttt{Gurobi} is a commercial solver that implements this BB method~\cite{gurobi}. 
This enumeration method provides rigorous global optimality and convergence guarantees, but the computational burden can grow exponentially with the problem size, limiting its applicability and scalability. 
Consequently, there is a substantial research interest focused on making the discrete search smaller and better guided through various algorithmic developments, such as logic-based search, cutting planes, and heuristics \cite{minlpreview2019}.
Among these efforts, Ising solvers are a promising approach for their potential of practical applicability to large-scale combinatorial optimization problems. 
It has been shown that many discrete optimization problems can be mapped to Ising models with only a polynomial overhead~\cite{cook_complexity_1971,karp_reducibility_1972}, and any advantage in solving these Ising problems over traditional methods would have significant practical implications \cite{Lucas2014,Mohseni_2022}.

Ising models, first developed to describe magnetism \cite{Lenz:460663,ising_beitrag_1925}, offer a physics-based representation of discrete optimization problems by encoding binary decision variables as spin states and defining the objective function as the energy of the system through a Hamiltonian function: 
\begin{equation}
  H = \sum_{i \in V(G)} h_i \sigma_i + \sum_{(ij)\in E(G)} J_{ij} \sigma_i \sigma_j
  \label{eq:ising_general}
\end{equation}
where \(\sigma_i \in \{-1, +1\}\) are binary spin variables, indexed by the vertices \(V(G)\) of graph \(G\), and the pairwise interactions are defined by the edges \(E(G)\) of the graph, where \(h_i\) and \(J_{ij}\) are the corresponding coefficients \cite{Berkley_2010,Bian2014}.
This formulation is mathematically equivalent to the quadratic unconstrained binary optimization (QUBO) model (Eq.~\eqref{eq:qubo_general}), and a simple variable transformation can be used to convert between the two representations (i.e. \(\sigma_i = 2x_i - 1\)).
\begin{equation}
\min_{x} \quad x^\top \mathbf{Q} x \quad \text{s.t.} \quad x \in \{0, 1\}^n
\label{eq:qubo_general}
\end{equation}
\noindent where \(\mathbf{Q}\) is without loss of generality an \(n\)-by-\(n\) square, symmetric matrix of coefficients, and \(x\) is a set of binary variables.
This equivalence establishes a connection between the Ising solvers and many constrained integer programming (IP) problems, as the IPs can be reformulated into QUBOs by incorporating the constraints into the objective function as quadratic penalty terms \cite{glover2019,Kochenberger2014}.
Ising solvers are specialized hardware that operate to find or approximate the ground states of the system, effectively providing solutions to the corresponding integer optimization problems \cite{Mohseni_2022,Lucas2014,toqubo:2023}. 
Many Ising hardware have been developed based on various technologies, including classical thermal annealers, quantum annealers, and dynamical system-based solvers such as the coherent Ising machines. 
A comprehensive review of these technologies and their implementations can be found in Reference \cite{Mohseni_2022}. 

In contrast to the deterministic classical solvers that strive to find and guarantee the global optimum, Ising solvers operate as stochastic heuristics, providing solution distributions rather than single optima.
Simulated annealing (SA) is a classical method that has been used to solve Ising or QUBO problems, exemplifying this approach by employing a probabilistic technique inspired by the physical annealing process. 
SA uses random sampling to search the solution space and makes probabilistic decisions to accept or reject new solutions based on a temperature parameter that gradually decreases over time, allowing the algorithm to escape local minima \cite{Henderson2003}.
The annealing techniques emulate the annealing process of metal thermal processing to attain the lowest lattice energy state; the analogy is that the optimal solution in a minimization optimization is found through a specific algorithmic treatment \cite{Kirkpatrick1983}.
While SA does not provide a global optimality guarantee, the randomized nature of search can lead to better exploration of the rugged energy landscape of the discrete optimization problems.
Building on this concept, various physics-based and physics-inspired solvers have emerged, providing other methods for solving Ising or QUBO problems. 
These systems encode the Hamiltonian directly into a physical system and exploit the system’s natural dynamics to drive it towards low-energy configurations.

Quantum annealers, such as D-Wave's Advantage systems, are specialized quantum computing devices designed to perform quantum annealing \cite{dwave_advantage}.
Quantum annealing (QA) directly implements Ising Hamiltonians as the native energy functions of the hardware to find the ground state of the encoded optimization problems through adiabatic quantum evolution \cite{Farhi2000,Rajak_2023}.
In this process, the system whose energy is described by an initial Hamiltonian is evolved into a final Hamiltonian slowly enough and without interactions with its environment, also known as \emph{adiabatically}.
Considering a time schedule $0 \leq s(t) \leq 1$, where $0 \leq t \leq t_f$ is the time, and the initial and final Hamiltonians, $H_0$ and $H_f$, respectively, a time-dependent Hamiltonian can be defined as:
\begin{equation}
    H(s) = A(s)H_0 + B(s)H_f ,
\end{equation}
\noindent where the functions $A$ and $B$ are smooth, differentiable, and satisfy the boundary conditions $A(1) = 0, B(0) = 0$.
In the simplest case, the time-dependent Hamiltonian can be a convex combination of the initial and final Hamiltonians in terms of the schedule, i.e., $A(s) = 1-s, B(s) = s$.
The quantum adiabatic approximation states that with a slow enough transition, the instantaneous energy ranking is tracked during the Hamiltonian evolution \cite{born1928beweis}.
If the initial energy state is the minimum according to the ranking, also known as the ground state, at the end of the evolution, the system will be sitting in the final Hamiltonian's ground state. 
The final Hamiltonian encodes the optimization problem of interest, and thus the ground state corresponds to the optimal solution of the minimization problem.
Other quantum algorithms, such as the quantum approximate optimization algorithm (QAOA), can also be implemented in the gate-based quantum computing framework to solve Ising or QUBO problems \cite{Farhi2014,Hadfield2017}. 
Interested readers can refer to Reference \cite{Sanders2020} for additional quantum algorithms for combinatorial optimization problems.

Another emerging paradigm is entropy computing (EC); it is an optic-based computing system introduced by Quantum Computing Incorporated (QCi). 
Dirac-1, originally presented as an optics-based platform for addressing Ising optimization problems~\cite{dirac1paper}, uses optical phases to represent the binary spin states and emulate the interactions between them through light interference and nonlinearity.
A feedback loop is employed to iteratively adjust the probability distribution towards the lower-energy states, with inherent noise playing a constructive role by promoting exploration and biasing the system toward favorable solutions rather than directly minimizing an energy function \cite{qci_eqc_2025}.
The emergence of these Ising solvers and other alternative computing hardware offers new opportunities to tackle combinatorial problems, and recent advances in this hardware have made sufficient progress to warrant exploring their potential in practical applications.

In process design and optimization, Ising solvers can be particularly useful for addressing the combinatorial aspects of design problems, such as equipment selection, process configuration, and operational scheduling. 
The primary focus of this work is on process superstructure optimization; however, the proposed framework and approach are not limited to this domain and are generalizable to other MINLP problems with decomposable structures. 
A process synthesis optimization problem can be written as a MINLP with an objective function subject to equality and inequality constraints, as shown below \cite{Mencarelli_2020,Turkay1996}.

\begin{equation}
\min_{x,y} f(x,y) \quad \text{s.t.} \quad h(x,y) = 0, \quad g(x,y) \leq 0, \quad x \in \mathbb{R}^n \subseteq X, \quad y \in \mathbb{Z}^m \subseteq Y 
\label{eq:mip}
\end{equation}

\noindent where $x$ represents continuous variables (e.g., flow rate), $y$ represents integer variables (e.g., unit selection). $f$ represents the objective function, which is the performance metric of the optimization problem (e.g., cost). 
The model equations $h$ describe the interaction of state variables with system physics (e.g., mass balances), and the inequalities $g$ describe specifications and operational or safety constraints (e.g., critical quality attribute requirements). 
In the case of maximization, the negative value of the corresponding subject of interest can be minimized (e.g., the negative value of the production mass). 
The set of constraints defines the feasible space of the problem, which is represented algebraically by $h$ and $g$.  
In cases where the derivation of such algebraic expressions is complex, simulation tools can be used to represent the process in a black-box optimization approach \cite{Casas-Orozco2021}.

Using a MINLP for end-to-end optimization (E2EO) has been actively explored in relatively well-behaved system applications, such as steady-state processes observed in the chemical industry \cite{Caballero2005, Navarro-Amoros2014, Corbetta2016}. 
However, in systems that require high modeling fidelity and exhibit complex behaviors, such as pharmaceutical processes, the mathematical programming approach to optimization has been limited in its applicability.
The derivation of equation-oriented optimization models for these systems is complex, and simulation-based optimization requires an iterative process that can result in a non-trivial computational burden \cite{BarhateRule2024}.
However, efforts have been made to integrate process synthesis optimization with overall process dynamics optimization. 
In Reference~\cite{BarhateRule2024}, the authors propose a rule-based and optimization-driven decision framework for optimizing flowsheets in a Drug Substance (DS) manufacturing process.
The methodology leverages heuristic rules, such as regulatory considerations and knowledge-based rules, as well as scenario analysis, to generate a smaller search space.
This framework efficiently narrows down the search space, thus evaluating alternative configurations more effectively overall.
However, the derivation of these rules is susceptible to user bias and lacks a quantitative evaluation of the alternative configurations.
In other words, this framework does not directly optimize the configuration selection itself. 

To address the complexity of discreteness and the desire for high fidelity in process design problems, we explore a heterogeneous computing approach that leverages the advantages of the Ising form and utilizes involved optimization algorithms, such as black-box optimization. 
Through decomposition, we isolate the discrete part of the problem, specifically the Ising part of the formulation, thereby offloading a source of hardness and better capturing the combinatorial aspect of the problem. 
Furthermore, we still maintain high fidelity and solve the nonlinear part of the problem using a black-box optimization approach. 
In this work, we investigate this approach and the application of the Ising solvers through both classical and quantum algorithms to solve combinatorial problems, with a focus on process synthesis and practical decision-making.

\subsection{Contributions of this work} 
The contributions of this work are as follows:
\begin{itemize}
  \item propose a framework for integrating Ising solvers into the solution algorithm of complex MINLP problems,
  \item application case studies of using Ising solvers, such as simulated annealing (SA), quantum annealing (QA), and entropy computing (EC), to solve discrete subproblems in process design problems,
  \item exploration of the use of novel computing methods, such as quantum annealing (QA) and entropy computing (EC), in process design problems,
  \item provide open-source code for the case studies to facilitate reproducibility and further research in this area\footnote{\url{https://github.com/SECQUOIA/pd_ising}}.
\end{itemize}

\section{Methods}
\label{sec:methods}
\begin{figure}
\centering
\includegraphics[width=12.5cm]{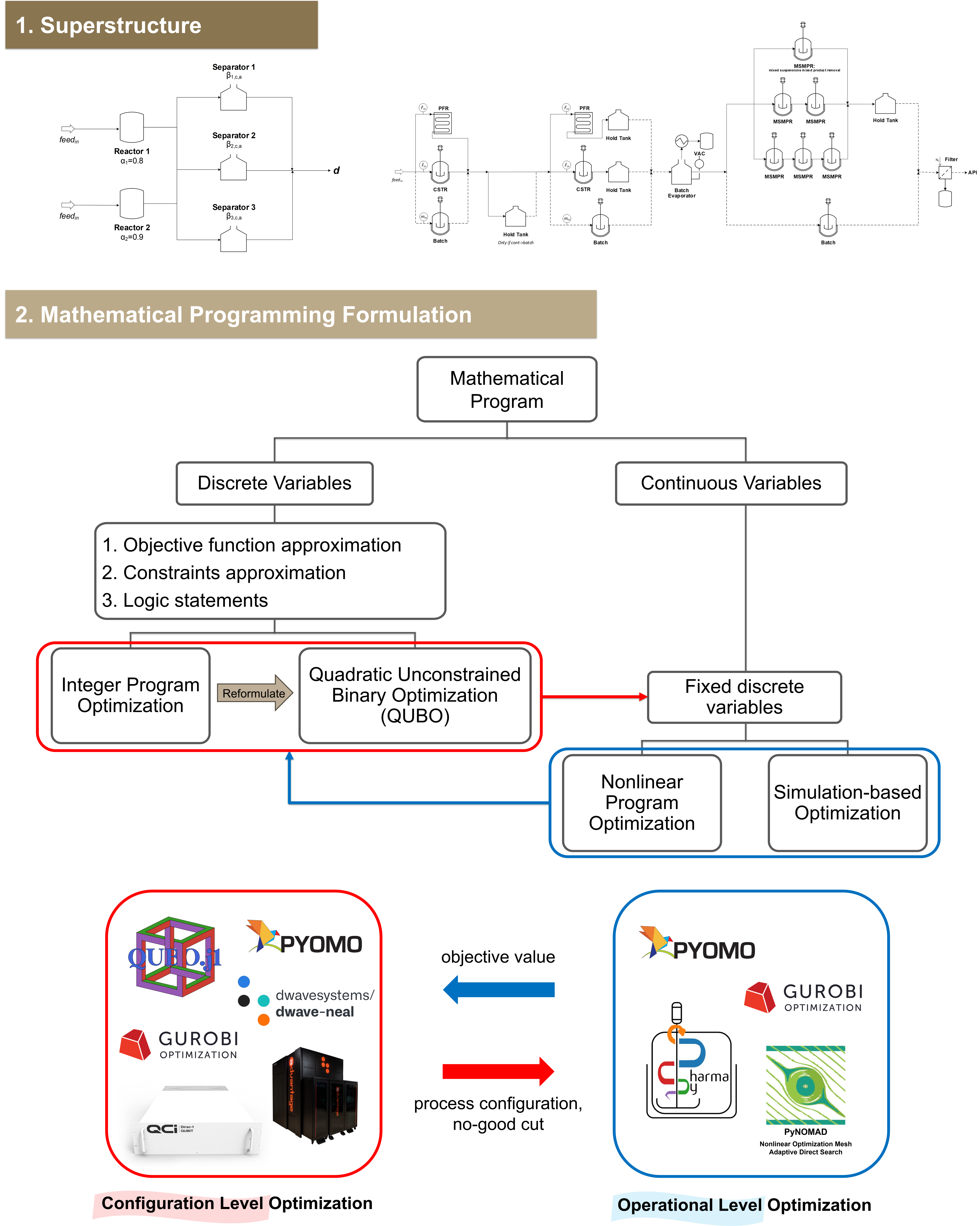} 
\caption{An illustrative workflow of the proposed approach. In step (1), we construct a single flowsheet of a superstructure that encapsulates the set of all feasible alternative process structures. Denoting critical flow paths and nodes, the process design and optimization is cast as a mathematical program, consisting of discrete (configuration) and continuous (operational) variables.  
Step (2), we decompose the problem into discrete and continuous parts, where the discrete subproblem is formulated as an integer program using approximation and logical statements. The integer program can further be reformulated into a quadratic unconstrained binary optimization to take advantage of Ising solvers. The solutions from the discrete subproblem fix the discrete variables in the continuous subproblem, where powerful nonlinear solvers and simulators can be leveraged. 
The bottom image illustrates the procedure used in the case studies presented in this work. The red box on the left represents discrete optimization at the configurational level, and the blue box on the right represents the continuous subproblem, which is handled by a classical solver or a simulator to optimize the process at the operational level.  
}
\label{fig:overall} 
\end{figure}

This section describes the methods used, including the formulation of discrete subproblems and the solution approaches employed. 
The continuous subproblems are not discussed in detail, as they are not the primary focus of this work; however, they are solved using a black-box optimization approach, either with a simulator or an optimization solver.
Two case studies are explored in this work: an illustrative example of ionic liquid selection and process configuration, and a more complex optimization problem for drug substance manufacturing. 
BB, SA, QA, and EC are used to solve the discrete subproblem in both case studies.
A workflow overview of the process is illustrated in Fig.~\ref{fig:overall}.

\subsection{Discrete Subproblem Formulation and Solution Methods}
As many discrete variables in decision-making MINLP problems are binary (e.g., unit selection), the discrete subproblem is formulated as an integer programming (IP) problem with binary variables in this work.
The formulation is implemented using \texttt{JuMP}\cite{jump} in \texttt{Julia} or using \texttt{Pyomo}\cite{bynum2021pyomo} in \texttt{Python} programming language and solved with \texttt{Gurobi}\cite{gurobi} for branch-and-bound method. 
Two modes of testing are implemented to find all feasible solutions. 
First, the problem is solved with the default solver settings (BB) that provide a single optimal solution.
To find all feasible solutions with this mode, an iterative procedure of adding a "no-good cut" is used to eliminate the solution of the current iteration.
In each iteration $n$, a "no-good cut" constraint is added to eliminate previous solution as infeasible as \( \sum_{y|y_{n-1} = 1} (1 - y) + \sum_{y|y_{n-1} = 0} y \geq 1 \quad \text{for } n = \{1, 2, \dots, n_{max}-1 \}\) with $n_{max} = N$ or the number of possible alternatives. 
Each iteration in this approach represents the solution of the combinatorial part of the MINLP, the fixing of the discrete variables, and the solution of the resulting NLP subproblem for evaluation. 
Secondly, the problem is solved using the solution pool feature of \texttt{Gurobi} (BB-Pool). 
This option prompts the solver to retain and return multiple solutions found along the search tree. 
The user can also specify the number of best solutions for the solver to return using this feature~\cite{gurobi}.
For these case studies, the total number of combinations of the binary variables is known, so this value is used as the setpoint for these parameters when operating in this mode.

For other Ising solvers, the IP subproblem is reformulated into a QUBO problem using the open-source \texttt{Julia} package \texttt{QUBO.jl} \cite{toqubo:2023}, and then solved via simulated annealing, quantum annealing, and entropy computing.
The simulated annealing (SA) method was implemented using the \texttt{dwave-neal} package \cite{dwave_neal} with default parameters, consisting of 1000 reads and 1000 sweeps. 
For quantum annealing (QA), the D-Wave Advantage 4.1 (QA-Adv1) and Advantage2 1.8 (QA-Adv2) quantum processing units (QPU)\cite{dwave_advantage} were used with the default annealing schedule and a sample anneal time of 20 $\mu$s, with 1000 reads \cite{Bian2014}.
The Advantage 4.1 QPU uses the Pegasus topology, while Advantage2 1.8 adopts the Zephyr topology, which provides higher native qubit connectivity and shorter minor-embedding chains for dense QUBO mappings.
We refer interested readers to References \cite{boothby_2016_chimera,Boothby2019_pegasus,Boothby2021_zephyr} for the details of the topologies and connectivity of the QPUs.
The embedding overhead in both case studies was examined to be less than 0.25 second and not considered separately in computational time results.
For entropy computing (EC), the Quantum Computing Inc. (QCi) Dirac-1 device was used with 100 samples per run, which was the maximum number allowed \cite{dirac1paper,dirac1patent}. 

For comparison in computation time, a performance metric called time-to-target (TTT) is introduced.
\begin{equation}
\text{TTT}_s = \tau \cdot \frac{\log(1 - s)}{\log(1 - p_{\text{target}})}
\end{equation}
\noindent where $\text{TTT}_s$ is the time required to achieve success with probability $s$ (typically 0.99) in reaching the \textit{target}.
$\tau$ is the execution time of the algorithm. $p_{\text{target}}$ is the probability of reaching the target solution.
If $p_{\text{target}} = 1$, which applies for deterministic solvers, then \(\text{TTT}_s = \tau\) \cite{Mcgeoch2013}.
Time to \textit{optimality} and finding \textit{all feasible solutions} are considered as targets in this work. 
Targeting optimality serves as a benchmark for the performance of the methods, while targeting all feasible solutions provides a measure of the solution diversity and completeness of the search space. 
From a decision-making perspective, finding all feasible solutions is relevant, if not crucial, as it allows for a comprehensive understanding of the alternatives and enables better-informed decisions.

The computational setup used in this work runs on a Linux Ubuntu 22.04 operating system and features an Intel\textsuperscript{\textregistered} i7-1365U processor with a base frequency of 1.80~GHz and 32.0~GB of RAM. 
The environment supports both Python~3.10.12 and Julia~1.11. Quantum computing resources include the D-Wave Advantage~4.1 quantum annealer and the QCi Dirac-1 entropy-based quantum computer. 
Optimization and modeling tasks were performed using the following packages: \texttt{JuMP}~v1.26.0, \texttt{Pyomo}~v6.7.3, \texttt{ToQUBO.jl}~v0.1.10, \texttt{PharmaPy}~v0.4.0, and \texttt{pyNOMAD}~v4.4.0. Solver and hardware interfaces include \texttt{Gurobi} (v11.0, v12.0.2), \texttt{dwave-neal}~v0.6.0, and \texttt{qci\_client}~v4.5.0.

\subsection{Quadratic Unconstrained Binary Optimization (QUBO) Reformulation}
A constrained IP problem with binary variables $y$ can be represented as Eq.~\eqref{eq:mip_std}; inequality constraints (e.g., $g$ in Eq.~\eqref{eq:mip}) can be written into \(Ay=b\) form by adding slack variables to the expressions \cite{glover2019}.
\begin{equation}
\min_{y} \quad c^\top y \quad \text{s.t.} \quad Ay = b \quad y \in \{0, 1\}^n 
\label{eq:mip_std}
\end{equation}
\noindent where \(c\) is the cost vector, \(A\) is the constraint matrix, and \(\mathbf{b}\) is the right-hand side vector.

These problems can be reformulated into a QUBO model by including quadratic infeasibility penalties in the objective function.
The constraints \(Ay=b\) are moved into the objective function with a penalty cost ($\rho$) as:
\begin{equation}
  \min_{y} \quad \mathbf{c}^T y + \rho(Ay - b)^\top(Ay - b)
  \label{eq:cx_mod}
\end{equation}
\noindent The penalty terms are: \(\rho(Ay - b)^\top(Ay - b) = \rho(y^\top(A^\top A)y - 2(A^\top b)y + b^\top b) \).
Taking advantage of \(y^2=y \text{ for } y \in \{0, 1\} \), the linear terms ($c$, $A^\top b$) appear on the diagonal of the matrix $Q$ in the general form of a QUBO as described in Eq.~\eqref{eq:qubo_general} \cite{glover2019}.

\subsection{Illustrative example: discrete subproblem formulation}\label{sec:ilprob}
The original problem is (P8) in Reference \cite{Iftakher}; its formulation is included in the Supplementary section for completeness. 
This case study presents an optimization problem involving the synthesis of a reactor-separator network with two reactors and three separator options while simultaneously selecting an ionic pair from a list of two cations and two anions for the process.
There are a total of 84 possible combinations of the discrete choices.
When formulated as a single monolithic problem, it is a MINLP with a nonlinear objective function and binary and continuous variables. 
This problem was decomposed into two subproblems: 1. the discrete network synthesis and ion pair selection problem, and 2. the continuous flow optimization at fixed discrete variables.
The discrete subproblem was formulated as an integer program (IP) with the objective function to minimize the total cost and binary variables for cation ($z_c$) and anion ($z_a$) selection, flow ($f_{i,j}$), and unit selection ($y_{r/s}$) as follows:

\begin{align}
\min_{y,w} \quad & \sum_{k \in \mathcal{K}} c^{\text{fixed}}_k y_k 
+ 2 \sum_{r \in \mathcal{R}} c^{\text{oper}}_r y_r \alpha_r + 2 \sum_{s \in \mathcal{S}} \sum_{c \in C} \sum_{a \in A} c^{\text{oper}}_s y_s \beta_{s,c,a} w_{c,a} \label{eq:ild_objf} \\[1em]
\text{s.t. } 
& f_{\text{src}, r} = y_r   \qquad \forall r \in \mathcal{R} \label{eq:ild_yr}\\
& f_{s, \text{sink}} = y_s  \qquad \forall s \in \mathcal{S} \label{eq:ild_ys}\\
& \sum_{r \in \mathcal{R}} f_{\text{src}, r} - f_{\text{src}, r_1} \cdot f_{\text{src}, r_2} = 1 \label{eq:ild_fsrc}\\
& \sum_{s \in \mathcal{S}} f_{s, \text{sink}} \geq 1 \label{eq:ild_fsink}\\
& f_{r,s} \cdot f_{\text{src}, r} = f_{r,s}, 
\quad f_{r,s} \cdot f_{s, \text{sink}} = f_{r,s} \qquad \forall r \in \mathcal{R}, s \in \mathcal{S} \label{eq:flow_conservation} \\
& (1 - f_{\text{src}, r}) + \sum_{s \in \mathcal{S}} f_{r, s} \geq 1 \qquad \forall r \in \mathcal{R} \label{eq:ild_rxtorflow}\\
& (1 - f_{s, \text{sink}}) + \sum_{r \in \mathcal{R}} f_{r,s} \geq 1 \qquad \forall s \in \mathcal{S} \label{eq:ild_sepflow}\\
& \sum_{c \in C} z_c = 1, \quad  \sum_{a \in A} z_a = 1 \label{eq:ild_onepair}\\
& w_{c,a} = z_c \cdot z_a \qquad \forall c \in C, a \in A \label{eq:ild_wrule}
\end{align}

\noindent where \(\mathcal{R}\) is the set of reactors, \(\mathcal{S}\) is the set of separators, and \(\mathcal{R}, \mathcal{S} \subset \mathcal{K}\) are the sets of units in the overall process. \(C\) is the set of cations, and \(A\) is the set of anions.
Additionally, $c^{\text{fixed}}_k$ is the fixed cost of unit $k$, $c^{\text{oper}}_{r/s}$ is the operating cost of reactor or separator $r/s$, $\alpha_r$ is the conversion factor of reactor $r$, and $\beta_{s,c,a}$ is the separation factor of separator $s$ for cation $c$ and anion $a$. 
Constraints for the subproblem are formulated with domain knowledge and logic-based statements; for example, flow conservation such as \emph{ if there is flow into a reactor, there must be flow into at least one of the separators} is enforced through a constraint, $ (1-f_{src,r})+\sum f_{r,s} \geq 1$, as shown in Eq. \eqref{eq:ild_rxtorflow}. 
Eqs. \eqref{eq:ild_yr} and \eqref{eq:ild_ys} ensure that the flow from the source to the reactor and from the separator to the sink is equal to the binary variable of unit selection, $y_r$ and $y_s$, respectively. 
Constraints \eqref{eq:ild_fsrc} and \eqref{eq:ild_fsink} ensure that there is at least one flow from the source to the reactor and from the separator to the sink, and constraints \eqref{eq:flow_conservation}-\eqref{eq:ild_sepflow} ensure that the flow conservation is satisfied at each unit. 
Lastly, a cation and an anion are selected using Eq. \eqref{eq:ild_onepair}, and the product of the two selections is represented by a binary variable $w_{c,a}$ in Eq. \eqref{eq:ild_wrule}.

\subsection{Simulation-based optimization problem example: discrete subproblem formulation}\label{sec:dsprob}
This case study presents a simulation-based optimization problem of a drug substance (DS) manufacturing process and serves as a case study of a more complex problem than the illustrative example.
The original problem is adopted from~\cite{Casas-Orozco2023} and describes a process for synthesizing a general drug substance, consisting of reaction, crystallization, and separation steps with multiple operating options for each step. 
The process configuration features two reactors, an evaporator for the solvent switch step in preparation for the crystallization step, and a filtration step to separate the solid active pharmaceutical ingredient (API). 
Three reactor types are considered for each reaction unit: a plug-flow reactor (PFR), and two continuously stirred tank reactors (CSTR) with different operating modes (continuous and batch). 
For evaporation, only the batch process is considered. Four options are considered for the crystallization step: one option is a batch crystallizer, and the other three options are continuous mixed suspension mixed product removal units (MSMPR) in series, ranging from one to three units. 
For each step, only one option is selected.

\begin{figure}
\centering
\includegraphics[width=15.5cm]{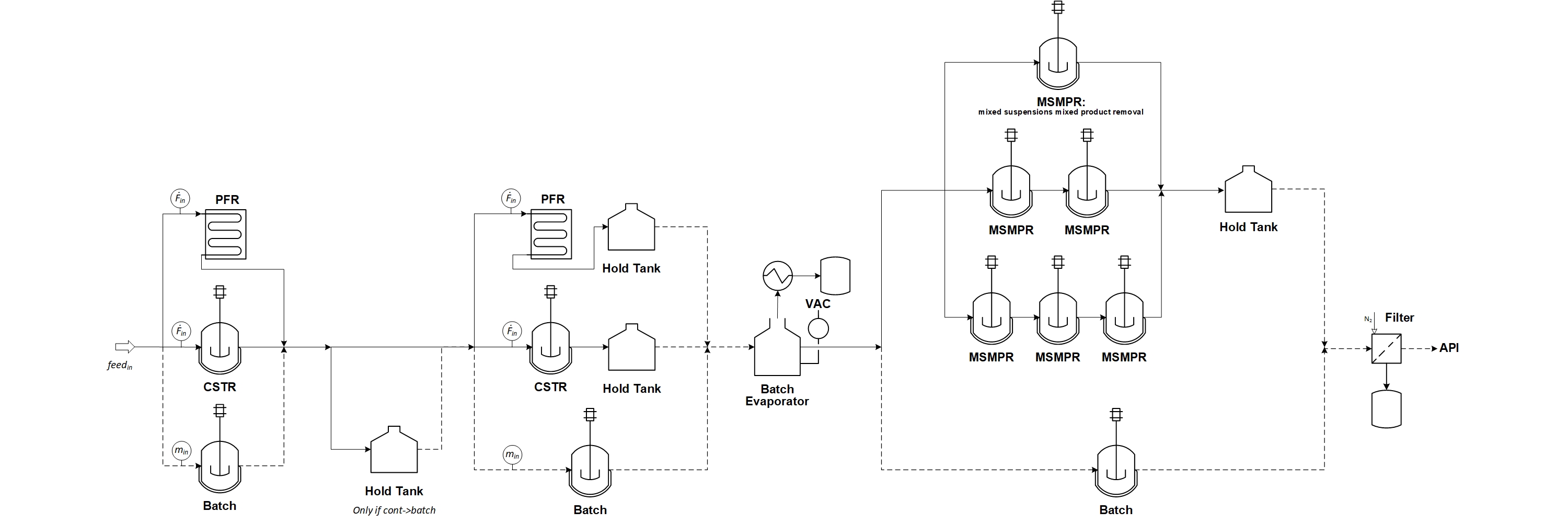} 
\includegraphics[width=15.5cm]{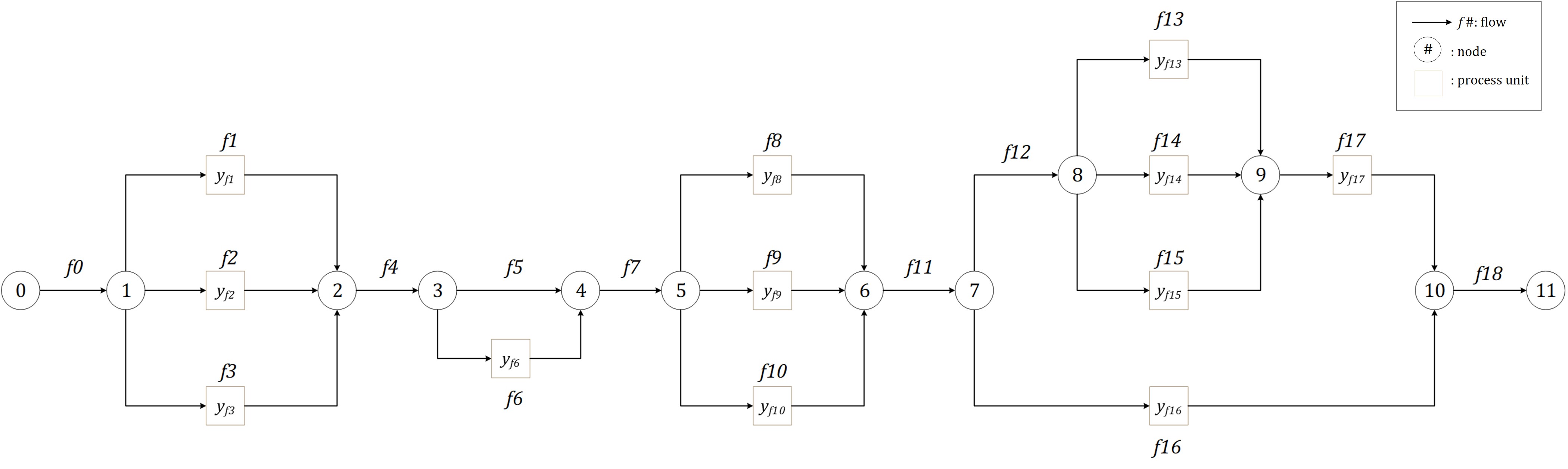} 
\caption{A superstructure of a drug substance manufacturing process (top) and a representative flow diagram for the main problem (bottom); dotted and solid lines indicate batch and continuous feed, respectively. }
\label{fig:ss} 
\end{figure}

The problem is decomposed into two subproblems: a discrete configuration selection problem and a simulation-based operational optimization problem. 
A superstructure of the process is shown in Fig.~\ref{fig:ss}. 
In total, there are 36 alternative configurations. 

The objective function of the configuration design subproblem is to minimize the total capital expense, and the objective function for the simulation-based subproblem is to maximize the production rate and crystal size. 
The simulation-based optimization framework, adopted from the literature \cite{BarhateHyb2024, BarhateRule2024, Laky2022}, is not discussed in detail; the primary focus of this work is on the discrete problem and the application of the heterogeneous computing framework.

Similarly to the illustrative example, the discrete problem is formulated as an integer program and as a QUBO problem, and solved using \texttt{Gurobi v.11.00} and the aforementioned heuristic methods (SA, QA, EC) and the corresponding parameters, respectively.
For the continuous subproblem, the configuration is fixed based on the solution of the discrete subproblem, and an open-package simulator, \texttt{PharmaPy}\cite{pharmapy}, is used to simulate the process and optimize the flow rates of the process, using \texttt{pyNOMAD}\cite{pynomad} as a black-box optimizer. 
The objective of the simulation-based framework is to maximize the production rate and crystal size. 

Two sets of variables are defined for optimization of the process configuration:
\(F = \{f_{00}, f_{01}, ...,f_{18} \}\), representing the binary flow variables, and \(L = \{y_{f_{01}}, y_{f_{02}},...,y_{f_{17}}\}\), representing the binary variables associated with the unit operation through which the corresponding flows pass, as illustrated in Fig.~\ref{fig:ss}.
For a system with parallel units or flows at disjunctions, these two variables ($f$ and $y$) may need to be treated separately.
In this case study, each disjunction represents a discrete choice of operating mode rather than a flow split; therefore, a single binary variable, $f$, is used to represent both flow and unit selection. 
The objective function of the subproblem is to minimize the total capital expense, $C$. 
The unit capital cost ($c_i$) can be estimated in various ways, such as using previously calculated data. In cases where data are unavailable, theory-based calculations or approximations suffice; methods such as the bare module method can be used to estimate the associated expense \cite{turton2018}.
For flows that do not activate any unit, the associated cost is set to zero. 

The formulation of the discrete subproblem is as follows:   

\begin{align}
\min_{f} \quad & \sum_{i \in F} c_i\, f_i  && \text{(objective)} \label{eq:objf}\\
\text{s.t.}\quad 
& \sum_{i \in \text{In}(n)} f_i \;=\; \sum_{i \in \text{Out}(n)} f_i && \forall\, n \in N \quad \text{(flow conservation)} \label{eq:mb}\\
% & g(f) \leq 0 \label{eq:inEC}  \\
& f_{source}=1,\quad f_{sink}=1 && \text{(source/sink activation)} \label{ss}\\
& g(f) \;\le\; 0 && \text{(continuous-to-batch constraints)} \label{ctb}\\
& f_i \in \{0,1\} && \forall\, i\in F. \label{bin}
\end{align}

\noindent where $F$ is the set of all flows, and $N$ is the set of nodes. 
Two types of constraints are declared: flow conservation at each node \eqref{eq:mb} and logic-based constraints \eqref{ctb}.
At each disjunction, only one selection of unit or flow rule is enforced through the mass conservation equality constraint.
Logic-based constraints, such as \emph{when a continuous unit is followed by a batch process, a holding tank must be placed in between these two units}, are incorporated into the optimization framework as inequality constraints (\(g(f) \leq 0\)).
For example, \( (1 - f_{01}) + (1 - f_{10}) + f_{06} \geq 1 \), this equation indicates that if the first reactor is PFR (\( f_{01} = 1 \)) and the second reactor is a batch reactor (\( f_{10} = 1 \)), then the holding tank must be installed (\( f_{06} = 1 \)).
The full formulation of the discrete subproblem is included in the Supplementary section (Section~\ref{sec:supp}) for completeness.

\section{Results}
\subsection{Illustrative example: Ionic Liquid Selection and Configuration of Reactor-Separator Network (IL)}
The results of the discrete subproblems are presented in this section, and the detailed formulation and approach are described in Section \ref{sec:ilprob}.

\begin{figure}
\centering
\includegraphics[width=16cm]{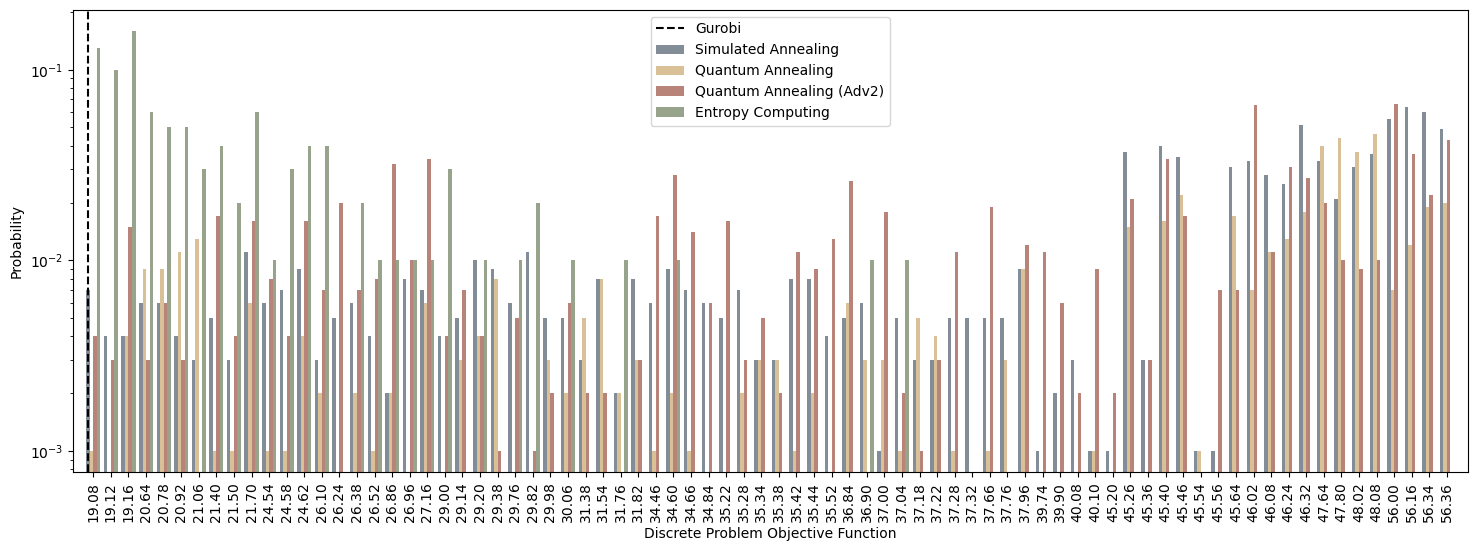} 
\caption{A plot of energy or probability of solutions found through various methods in the IL case study. The vertical line indicates the optimal solution of the discrete subproblem found by \texttt{Gurobi}. Infeasible solutions are not shown in this plot.}
\label{fig:energyplot} 
\end{figure}

The bar graph results presented in Fig.~\ref{fig:energyplot} show the energy or probability of solutions found in a single run of the heuristic methods on the y-axis, ranked by the discrete problem's objective function value on the x-axis.  
The plot illustrates that out of 84 possible configurations identified by Gurobi, SA was able to find all 84 feasible solutions in a single run. 
In contrast, QA-Adv1, QA-Adv2 and EC found 66, 78, and 30 feasible solutions, respectively. 
All methods were able to find the optimal solution.
The algorithm execution time values along with the calculated time-to-target (TTT) metrics are summarized in Table~\ref{tab:il-tts}. 
Total time for QA-Adv2 was not included in the analysis due to a post-processing error that led to incorrect values. 
Overall, BB achieved the fastest computation time for both finding the optimal solution and identifying all feasible solutions, followed by SA, QA, and EC. 
The solution pool mode of BB method (BB-Pool) resulted in the longest time to find all feasible solutions. 

These results highlight the differences among these methods. The heuristic methods can explore all feasible solutions in a single run, but require longer times to do so. 
In contrast, BB quickly solves for optimality, albeit with multiple iterations corresponding to the problem size, to discover all feasible solutions.
When BB method is applied to search for \textit{all} or multiple feasible solutions in BB-Pool mode, it becomes orders of magnitude slower than the heuristic methods.
Compared to the default mode with iterations, BB-Pool does not have access to information such as no-good cuts to navigate the search direction, slowing down the discovery of feasible solutions. 
Furthermore, the introduction of a new binary variable, $f$, in formulation of the discrete problem further enlarges the state-space. 
Since pruning branches is only possible through infeasibility in this mode, enumeration over the resulting large search space becomes inefficient. 

The optimal discrete choices from the original problem were found in the 72nd iteration solution of Gurobi for the discrete subproblem.
This mismatch occurred from the misalignment of the objective functions in the two problems. 
Although the objective function of the discrete subproblem was formulated to approximate that of the original problem, it could not account for the complexities involved with the continuous variables in the original problem. 
As a result, the subproblem solutions did not replicate the exact ranking of solutions in the iterations. 
Nevertheless, the integer program formulation of the discrete subproblem allowed the use of the discrete solvers and other optimization techniques.

\subsection{Simulation-based optimization problem example: Drug Substance Manufacturing Process Optimization (DSMFG)}
\begin{figure}
\centering
\includegraphics[width=14cm]{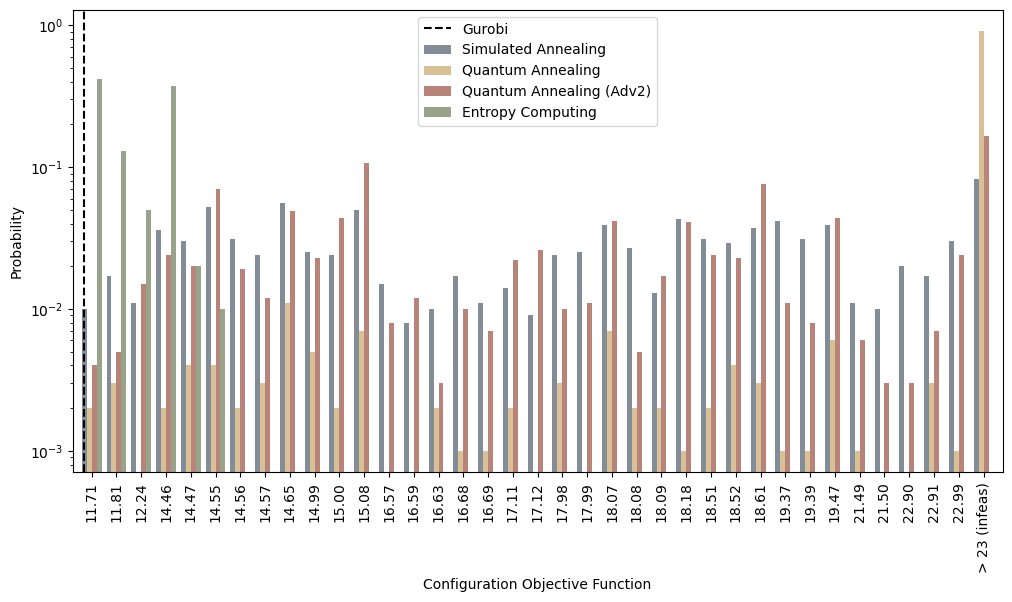} 
\caption{A plot of energy or probability of feasible solutions found through annealing methods. The vertical line indicates the optimal solution of the discrete subproblem, indicated by the first result of \texttt{Gurobi}. Infeasible solutions are not shown in this plot.}
\label{fig:energyplot_sim} 
\end{figure}

\begin{figure}[htbp]
\centering
\includegraphics[width=11cm]{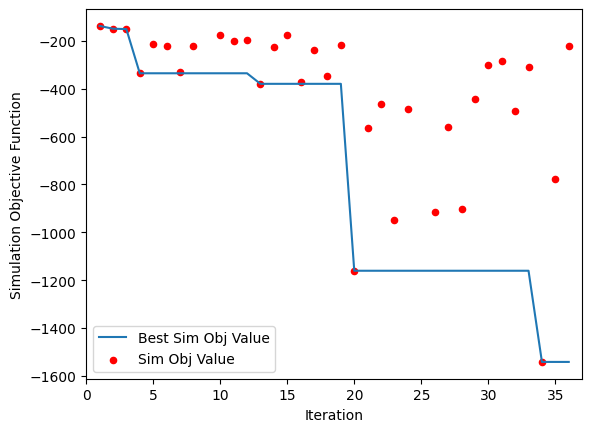} 
\caption{Simulation objective function values ranked according to discrete subproblem solution by IP iterations (x-axis). Red dots indicate the objective function value of the corresponding simulation, and the blue line shows the best objective function value found with each iteration.}
\label{fig:ds-sim} 
\vspace{-8pt}
\end{figure}

All results of the discrete subproblem are plotted in Fig.~\ref{fig:energyplot_sim}; the simulation objective values at each iteration of IP are plotted along with the best objective value found with the iterations in Fig.~\ref{fig:ds-sim}. The TTT metrics are presented in Table~\ref{tab:ds-tts}.

Similar trends were observed in this case study as the illustrative example.
All methods were able to find the optimal solution. 
With a single run, SA and QA-Adv2 identified all 36 feasible solutions, while QA-Adv1 and EC found 28 and 6 feasible solutions, respectively. 
BB in both default (BB) and solution pool (BB-Pool) modes achieved the fastest computation time for both finding the optimal solution and identifying all feasible solutions, while entropy computing exhibited the slowest computation time. 
Unlike in the illustrative example, the BB-pool method outperformed the heuristic methods in computational speed in this case study. 
This result is understandable given that the discrete subproblem in this case study is smaller, leading to a smaller search space. 

\begin{figure}[htbp]
\centering
\includegraphics[width=12.5cm]{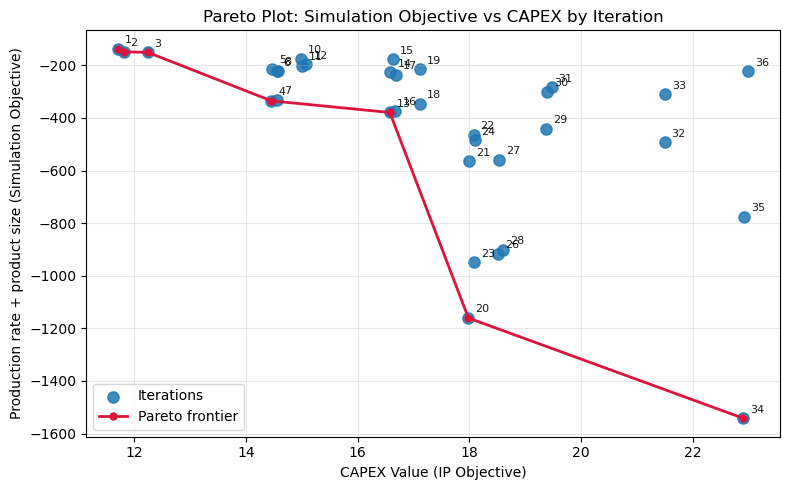}
\caption{Pareto plot of simulation objective function values (y-axis) versus IP discrete subproblem objective values (x-axis) for the DSMFG case study; in-plot labels indicate the corresponding iteration numbers. The simulation objective maximizes the production rate and crystal size of the final product, with more negative values indicating better performance, while the discrete subproblem objective minimizes the capital cost of the process configuration. The Pareto front (red) highlights the trade-off between productivity and product quality and cost.}
\label{fig:paretoplot}
\end{figure}

The best simulation objective value was achieved at the 34th iteration in the BB method. 
The corresponding configuration for the process was to choose the plug-flow reactor, continuously stirred tank reactor, and batch process for the first reaction, second reaction, and the crystallization unit selections respectively. 
The optimal configuration from the discrete subproblem was batch reactors for both reactions followed by a single continuous crystallization unit.
Because the discrete subproblem minimizes capital cost while the simulation objective maximizes production and product quality, the opposing rankings between the two subproblems are consistent with the expected trade-off between cost and performance.
Multi-objective optimization can simultaneously optimize the overall system for several purposes (e.g., minimizing capital cost while maximizing product purity and production rate), and pareto plots like Figure~\ref{fig:paretoplot} can be used to visualize the overall trade-off between the two objectives to better guide decision-making when selecting an optimal solution. 

\section{Discussion}
In general performance, the heuristic methods (SA, QA, EC) are slower than BB in finding the optimal solution.
Among the novel computing approaches, EC showed slower computation time and identified fewer solutions compared to QA. 
While QA identified solutions that were more broadly distributed across the solution space, EC solutions were more clustered near the optimal solution, indicating that EC tends to find the high-quality solutions early, but explores the solution space less extensively. 
These observations highlight a fundamental trade-off between solution optimality and diversity among these methods. 
Heuristic approaches can explore the solution space more broadly and provide diverse solutions within a single run. 
In contrast, BB requires iterations and additional mechanisms like cuts to recover multiple feasible solutions, and as the problem scale grows, the computational effort required for such enumeration increases exponentially. 
The results from the BB-Pool method further demonstrate how the growth of the search space increases the computational burden of the BB method.

A clear understanding of these methodological differences can enable their more effective use in practical applications.
For practical decision-making processes, the size of the optimization problem is often large, and exploring all combinatorial solutions is intractable and inefficient. 
In cases like this, heuristic methods can be employed to quickly explore the solution space and identify feasible discrete solutions. 
Several decisions can be made based on the heuristics and practical constraints, such as cost, time, and resources. 
Then the remaining continuous problem can be solved to complete the optimization process. 
For instance, in a process design problem, with a single run of heuristics, all the configuration decisions can be identified, and impractical process choices can be eliminated with user's domain knowledge in the first step. 
Then, the selected few configurations may proceed with simulations in parallel to gain a comprehensive understanding of the process implications without loss of fidelity. 

Selecting the appropriate computing method for specific optimization problems can significantly impact the computational effort required to solve them. 
In the IL case study, several degenerate solutions were observed in the subproblem, where multiple distinct feasible solutions resulted in the same objective function value. 
This illustrates one of the key differences between heuristic methods and deterministic solvers, such as \texttt{Gurobi}. 
In cases of degeneracy, deterministic solvers cannot prune any branches with the same objective value, and in the worst-case scenario, would have to evaluate all combinations of the discrete choices, which can increase exponentially with the problem size.  
Even worse, if cuts based on the objective value were imposed instead of an integer cut, there is a risk of eliminating alternative feasible, even optimal, solutions to the original problem.  
Additionally, for a simple problem such as the IL selection, current classical solvers, such as \texttt{Gurobi v.12.0}, can solve the original MINLP without modifications to the nonlinear objective function.
Further processing of the problem into a QUBO problem actually introduces additional slack variables and constraints, unnecessarily increasing the problem size and complexity. 
In contrast, complex problems such as drug substance manufacturing process optimization, characterized by differential-algebraic equations, nonlinearities, and discrete decisions, are intractable when approached as a single monolithic MINLP. 
By applying a decomposition strategy, we can systematically integrate discrete configuration decisions into large-scale problems and find practical solutions. 

Lastly, it is essential to note that since these are different methods, different parameters were used, and their performance results may not be directly comparable. 
For heuristic methods, the parameters were set to their default values, and a limit of 100 samples was imposed for EC. 
There is a potential for performance improvement through parameter tuning and by implementing an iterative process of adding cuts similar to that used in the BB method \cite{Bernal2024-benchmarking}.

\section{Conclusion}
This work proposes a heterogeneous computing approach for solving MINLPs through a decomposition strategy with Ising solvers and other NLP solution methods. 
The IP formulation was first used to represent the discrete subproblem for BB, and the QUBO reformulation was applied to express the discrete problem in a form suitable for Ising solvers.
Classical and Ising-based solvers were explored and compared for their computational performance and solution quality in two case studies: one simple and one complex, of process design optimization.
Lastly, the open-source code for the case studies are provided for general use and application of the framework and promote further research in the area. 

The differences in algorithm execution among these methods, along with the optimality guarantee and solution diversity of their results, were discussed in the context of practical decision-making purpose. 
The classical commercial solver enabled us to find the global optimal solution to the subproblem, yet it required iterations to identify additional feasible solutions. 
Other heuristic methods provided a probability distribution of solutions with a single execution of the algorithms, but not all methods were able to find the entire set of feasible solutions. 
Simulated annealing (SA) was able to find all feasible solutions, while the novel computing methods (QA and EC) did not find all feasible solutions in a single run. 
These results highlight that a better understanding of the differences in solution methods can help in more effectively leveraging these methods and hardware in practice. 
As a last remark, even though these emerging technologies do not match the classical solvers in convergence speed, their complementary strengths in scalability and solution exploration, along with the continuing developments in the hardware and algorithms, warrant continued investigation in future research. 

\section*{Acknowledgments}
This material is based upon work supported by the Center for Quantum Technologies under the Industry-University Cooperative Research Center Program at the US National Science Foundation under Grant No. 2224960.
We also acknowledge the financial support of Quantum Computing Inc. in sponsoring this research.

We would also like to acknowledge Daniel Laky and Daniel Casas-Orozco for their support and help in using PharmaPy and associated code. 

\clearpage
\section*{Tables}

\begin{table}[H]
\centering
\caption{Summary table of time to target (Opt: optimality, Feas: all feasible solutions) results of BB, SA, QA (D-Wave Advantage, Advantage2), and EC (QCi Dirac 1) for the discrete subproblem of the illustrative problem. All times are in seconds. The subscript "quantum" indicates the time required for the quantum processing unit only, while "total" includes both communication overhead and the time for the quantum processing unit. The subscript "device" indicates the overall computing time\cite{qci_dirac3}.}
\label{tab:il-tts}
\renewcommand{\arraystretch}{1.5}
\begin{tabular}{@{}>{\centering\arraybackslash}p{5.2cm}
                >{\centering\arraybackslash}p{2.7cm}
                >{\centering\arraybackslash}p{2.7cm}
                >{\centering\arraybackslash}p{2.7cm}@{}}
\toprule
& \textbf{Execution Time} & \multicolumn{2}{c}{\textbf{Time to Target ($TTT_{99}$)}} \\
\cmidrule(lr){2-4}
\textbf{Solution Method} & \textbf{$\tau$} & \textbf{$TTOpt_{99}$} & \textbf{$TTFeas_{99}$} \\
\midrule
\textbf{BB (Gurobi)}            & 0.003  & 0.003 ($=\tau$)  & 0.477 \\
\textbf{BB-Pool (Gurobi Pool)} & 62.95 & 0.003 & 62.95 \\
\textbf{SA}                     & 0.34   & 222.9            & 0.54 \\
\textbf{QA-Adv1$_{\text{quantum}}$}  & 0.13   & 620.5            & - \\
\textbf{QA-Adv1$_{\text{total}}$}    & 1.59   & 7313.9           & - \\
\textbf{QA-Adv2$_{\text{quantum}}$}  & 0.16   & 188.5           & - \\
\textbf{EC$_{\text{device}}$} & 37.0     & 1223.5           & - \\
\bottomrule
\end{tabular}
\end{table}

\begin{table}[H]
\centering
\caption{Summary table of time to target (Opt: optimality, Feas: all feasible solutions) results of BB, SA, QA (D-Wave Advantage, Advantage2), and EC (QCi Dirac 1) for the discrete subproblem of the drug substance manufacturing problem. All times are in seconds. The subscript "quantum" indicates the time required for the quantum processing unit only, while "total" includes both communication overhead and the time for the quantum processing unit. The subscript "device" indicates the overall computing time\cite{qci_dirac3}.}
\label{tab:ds-tts}
\renewcommand{\arraystretch}{1.5}
\begin{tabular}{@{}>{\centering\arraybackslash}p{5.2cm}
                >{\centering\arraybackslash}p{2.7cm}
                >{\centering\arraybackslash}p{2.7cm}
                >{\centering\arraybackslash}p{2.7cm}@{}}
\toprule
& \textbf{Execution Time} & \multicolumn{2}{c}{\textbf{Time to Target ($TTT_{99}$)}} \\
\cmidrule(lr){2-4}
\textbf{Solution Method} & \textbf{$\tau$} & \textbf{$TTOpt_{99}$} & \textbf{$TTFeas_{99}$} \\
\midrule
\textbf{BB (Gurobi)} & 0.0009 & 0.0009 ($=\tau$) & 0.0790 \\
\textbf{BB-Pool (Gurobi Pool)} & 0.073 & 0.0009 & 0.073 \\
\textbf{SA} & 0.606 & 197.8 & 1.1 \\
\textbf{QA-Adv1$_{\text{quantum}}$} & 0.136 & 311.9 & - \\
\textbf{QA-Adv1$_{\text{total}}$} & 0.410 & 942.7 & - \\
\textbf{QA-Adv2$_{\text{quantum}}$} & 0.174   & 199.6 & 0.44 \\
\textbf{EC$_{\text{device}}$} & 35.0 & 295.9 & - \\
\bottomrule
\end{tabular}
\end{table}

\clearpage

\bibliographystyle{unsrt}
\bibliography{bibliography}
\clearpage

\section{Supplementary}\label{sec:supp}
\subsection{Ionic Liquid Selection and Configuration of Reactor Separator Network Optimization Problem (original MINLP formulation)}

\begin{align}
  \min \quad & 
  \sum_{k \in K} c_k^f y_k 
  + \sum_{i \in I_{K_r}^{\text{in}}} c_i^I x_i^{0.6} 
  + \sum_{k \in K_s} c_k^I ( \sum_{i \in I_{K_s}^{\text{in}}} x_i )^2 
  + \sum_{k \in K_s} c_k^e ( \sum_{i \in I_{K_s}^{\text{in}}} x_i - \sum_{i \in I_{K_s}^{\text{out}}} x_i )
  \label{eq:il_objf} \\
  \text{s.t.} \quad & f_k^L y_k \leq \sum_{i \in I_k^{\text{in}}} x_i \leq f_k^U y_k, \qquad \forall k \in K \label{eq:il_fbound} \\
  & \sum_{i \in I_{K_r}^{\text{out}}} x_i = \alpha_k \sum_{i \in I_{K_r}^{\text{in}}} x_i, \qquad \forall k \in K_r \label{eq:il_rxt} \\
  & \sum_{c \in \text{Cat}} z_c = 1, \qquad \sum_{a \in \text{An}} z_a = 1 \label{eq:il-oneion} \\
  & x_i \geq \beta_{k, c, a} \sum_{i \in I_{K_s}^{\text{in}}} x_i - M (2 - z_c - z_a), \qquad 
  \forall i \in I_{K_s}^{\text{out}},\, k \in K_s,\, c \in \text{Cat},\, a \in \text{An} \label{eq:il-sep1} \\
  & x_i \leq \beta_{k, c, a} \sum_{i \in I_{K_s}^{\text{in}}} x_i + M (2 - z_c - z_a), \qquad 
  \forall i \in I_{K_s}^{\text{out}},\, k \in K_s,\, c \in \text{Cat},\, a \in \text{An} \label{eq:il-sep2} \\
  & \sum_{i \in I_{K_s}^{\text{out}}} x_i \geq d \label{eq:il-demand}
\end{align}

\begin{align*}
\text{Cat} &= \{ c_1, c_2, \dots, c_C \} & \text{set of cations} \\
\text{An} &= \{ a_1, a_2, \dots, a_A \} & \text{set of anions} \\
\text{IL} &= \{ c \times a \mid c \in \text{Cat},\, a \in \text{An} \} & \text{set of all feasible ionic liquids (ILs)} \\
K &= \text{set of all process units} \\
K_r &\subset K & \text{set of all reactors} \\
K_s &\subset K & \text{set of all separators} \\
I_k^{\text{in}} &= \text{set of all inlet streams to unit } k \\
I_k^{\text{out}} &= \text{set of all outlet streams from unit } k \\
d &= \text{demand requirement} \\
c_k^f &\quad \text{fixed cost for unit } k \\
c_i^I &\quad \text{operating cost for reactor stream } i \in I_{K_r}^{\text{in}} \\
c_k^I &\quad \text{operating cost for separator } k \\
c_k^e &\quad \text{emissions or waste cost for separator } k    \\
\alpha_k &\quad \text{stoichiometric reactor conversion factor for reactor } k \in K_r \\
\beta_{k,c,a} &\quad \text{separation factor for separator } k \in K_s \text{ with IL } (c,a) \\
x_i &\quad \text{flowrate of stream } i \\
z_c &\in \{0,1\} \quad \text{binary: cation } c \text{ selected} \\
z_a &\in \{0,1\} \quad \text{binary: anion } a \text{ selected} \\
y_k &\in \{0,1\} \quad \text{binary: unit } k \text{ selected}
\end{align*}
\newpage

\subsection{Drug Substance Manufacturing Process Optimization Problem (discrete subproblem IP formulation)}

\begin{align}
\min_{f} \quad & \sum_{i \in F} c_i\, f_i  && \text{(objective)} \label{eq:objf_ds}\\
\text{s.t.}\quad 
& \sum_{i \in \text{In}(n)} f_i \;=\; \sum_{j \in \text{Out}(n)} f_j && \forall\, n \in N \quad \text{(flow conservation)} \label{eq:mb_ds}\\
& (1-f_{01})+(1-f_{10})+f_{06} \;\ge\; 1 \label{ctb1_ds}\\
& (1-f_{02})+(1-f_{10})+f_{06} \;\ge\; 1 \label{ctb2_ds}\\
& (1-f_{03})+(1-f_{10})+(1-f_{06}) \;\ge\; 1 \label{ctb3_ds}\\
& f_{01}+f_{02}+f_{05} \;\ge\; 1 \label{ctb4_ds}\\
& f_{05}+f_{10} \;\ge\; 1 \label{ctb5_ds}\\
& f_{00}=1,\quad f_{18}=1 && \text{(source/sink activation)} \label{ss_ds}\\
& f_i \in \{0,1\} && \forall\, i\in F. \label{bin_ds}
\end{align}

\[
\begin{aligned}
F &:\ \text{set of flows}, \\
N &:\ \text{set of nodes},\\
\mathrm{In}(n) &\subseteq F:\ \text{incoming flows to node } n\in N, \\
\mathrm{Out}(n) &\subseteq F:\ \text{outgoing flows from node } n\in N.
\end{aligned}
\]

\[
\begin{aligned}
c_i &\ge 0 &&\text{cost of flow } i,\quad i\in F,\\
f_i &\in \{0,1\} &&\text{binary decision variable (1 if flow $i$ is selected)},\quad i\in F.
\end{aligned}
\]

\end{document}